\documentclass[11pt,reqno]{amsart}
\usepackage{tikz}

\usepackage{color}
\usepackage{hyperref}
\hypersetup{colorlinks=true,linkcolor=blue, citecolor=blue}
\usepackage{curves}
\usepackage{amsmath,amsthm,amscd,amsfonts,amssymb,euscript}
\usepackage{graphics}
\usepackage{pdfpages}
\usepackage{mathtools}

\usepackage{latexsym}
\usepackage{color} 
\theoremstyle{plain}

\font\manual=manfnt
\def\dbend{{\manual\char127}} 

\def\danger{\begin{trivlist}\begin{footnotesize}\item[]\noindent%
\begingroup\hangindent=3pc\hangafter=-2
\def\par{\endgraf\endgroup}%
\hbox to0pt{\hskip-\hangindent\dbend\hfill}\ignorespaces}
\def\enddanger{\par\end{footnotesize}\end{trivlist}}

\def\ddanger{\begin{trivlist}\begin{footnotesize}\item[]\noindent%
\begingroup\hangindent=3pc\hangafter=-2
\def\par{\endgraf\endgroup}%
\hbox to0pt{\hskip-\hangindent\dbend\kern2pt\dbend\hfill}\ignorespaces}
\def\endddanger{\par\end{footnotesize}\end{trivlist}}

\DeclareFontFamily{OT1}{rsfs}{}
\DeclareFontShape{OT1}{rsfs}{n}{it}{<-> rsfs10}{}
\DeclareMathAlphabet{\mathscr}{OT1}{rsfs}{n}{it}

\begin{document}

\numberwithin{equation}{section}

\newtheorem{guess}{\sc Theorem}[section]
\newcommand{\bth}{\begin{guess}$\!\!\!${\bf }~}
\newcommand{\eeth}{\end{guess}}

\newtheorem{propo}[guess]{\sc Proposition}

\newcommand{\bprop}{\begin{propo}$\!\!\!${\bf }~}
\newcommand{\eprop}{\end{propo}}

\newtheorem{lema}[guess]{\sc Lemma}
\newcommand{\blem}{\begin{lema}$\!\!\!${\bf }~}
\newcommand{\elem}{\end{lema}}

\newtheorem{defe}[guess]{\sc Definition}
\newcommand{\bdefe}{\begin{defe}$\!\!\!${\it }~}
\newcommand{\edefe}{\end{defe}}

\newtheorem{coro}[guess]{\sc Corollary}
\newcommand{\bcor}{\begin{coro}$\!\!\!${\bf }~}
\newcommand{\ecor}{\end{coro}}

\newtheorem{rema}[guess]{\it Remark}
\newcommand{\brem}{\begin{rema}$\!\!\!${\it }~\rm}
\newcommand{\erem}{\end{rema}}

\theoremstyle{remark}
\newtheorem{assump}[guess]{Assumption}

\newcommand{\spec}{{\rm Spec}\,}
\newtheorem{notation}{Notation}[section]
\newcommand{\bnot}{\begin{notation}$\!\!\!${\bf }~~\rm}
\newcommand{\enot}{\end{notation}}

\newcommand{\bpr}{\begin{proof}}
\newcommand{\epr}{\end{proof}}

\numberwithin{equation}{guess} 
\newcommand{\beqa}{\begin{eqnarray}}
\newcommand{\eeqa}{\end{eqnarray}}
\newtheorem{thm}{Theorem}
\theoremstyle{definition}
\newtheorem{say}[guess]{\bf}
\newtheorem{hint}[thm]{Hint}
\newtheorem{example}[thm]{Example}
\newcommand{\bsem}{\begin{say}$\!\!\!${\it }~~\rm}
\newcommand{\esem}{\end{say}}
\newtheorem{observe}[subsubsection]{Observation}

\newcommand{\ha}{\sf h}
\newcommand{\g}{\sf g}
\newcommand{\ta}{\sf t}
\newcommand{\s}{\sf s}
\newcommand{\ctext}[1]{\makebox(0,0){#1}}
\setlength{\unitlength}{0.1mm}

\newcommand{\wt}{\widetilde}
\newcommand{\Lr}{\Longrightarrow}
\newcommand{\Aut}{\mbox{{\rm Aut}$\,$}}
\newcommand{\ul}{\underline}
\newcommand{\ol}{\bar}
\newcommand{\lr}{\longrightarrow}
\newcommand{\sh}{{\sf h}}

\newcommand{\ba}{{\mathbb A}}

\newcommand{\bc}{{\mathbb C}}
\newcommand{\bp}{{\mathbb P}}
\newcommand{\bz}{{\mathbb Z}}
\newcommand{\bq}{{\mathbb Q}}
\newcommand{\bn}{{\mathbb N}}
\newcommand{\bg}{{\mathbb G}}
\newcommand{\br}{{\mathbb R}}
\newcommand{{\bh}}{{\mathbb H}}
\newcommand{\bo}{{\bar \omega}'}
\newcommand{\po}{{\omega}'}

\newcommand{\ct}{{\mathcal T}}
\newcommand{\cc}{{\mathcal C}}
\newcommand{\cl}{{\mathcal L}}
\newcommand{\cv}{{\mathcal V}}
\newcommand{\cf}{{\mathcal F}}
\newcommand{\cb}{{\Lambda}}
\newcommand{\ch}{{\mathcal H}}

\newcommand{\mfc}{{\sf C}}
\newcommand{\ce}{{\mathcal E}}
\newcommand{\co}{{\mathcal O}}

\newcommand{\cm}{{\mathcal M}}

\newcommand{\cs}{{\mathcal S}}
\newcommand{\cg}{{\mathcal G}}
\newcommand{\ca}{{\mathcal A}}
\newcommand{\hra}{\hookrightarrow}
\newcommand{\mfu}{{\sf U}}

\newtheorem{ack}{\it Acknowledgments}       
\renewcommand{\theack}{} 

\title{Parabolic Bundles and Parabolic Higgs bundles}

\author{D.S. Nagaraj} 
\address{Institute of Mathematical Sciences, Taramani, Chennai-600115, India, dsn@imsc.res.in}

\begin{abstract}
This is a survey article about parabolic bundles and parabolic Higgs bundles.
\end{abstract}
\maketitle

\vspace{2mm}
\noindent

\section{Introduction}

The main aim of this article is to give a sample survey about Parabolic 
bundles and Parabolic Higgs bundles. This is an expanded version of my 
survey  talk 
at the conference NS@50 held at Chennai Mathematical Institute during
the month of October 2015. There will be many glaring omissions which is 
partially due time constraint and mainly due to my limited knowledge 
of the subject.

\section{Parabolic vector bundles on a Riemann Surface }

\bdefe\label{def1} Let $X$ be smooth irreducible curve over an 
algebraically closed field $K.$ Let  $p_1,\ldots , p_n$ be fixed
finite set of closed points of $X.$  A parabolic vector bundle on $X$ 
with parabolic structure at 
$p_1,\ldots , p_n$ is a vector bundle $W$  on $X$ together with
the following data at each $p=p_i,$ 
\begin{itemize}
\item[a)] a flag 
	$$W_p =F_1W_p\supset F_2W_p\supset \ldots \supset F_rW_p,$$
of subspaces of the vector space $W_p,$ the fiber of $W$ at $p$
\item[b)] and real weights 
 $$\alpha_1,\ldots, \alpha_r$$
  attached to $F_1W_p,\ldots, F_rW_p$  such that
       $0\leq \alpha_1<\ldots < \alpha_r<1.$	
\end{itemize}	
\edefe
   
\brem
0) By abuse of notation we say that $W$ is a parabolic vector bundle
with parabolic structure at $p_1,\ldots , p_n.$

1) The numbers $k_1=\text{dim}F_1W_p-\text{dim}F_2W_p,\dots,
k_r=\text{dim}F_rW_p$ are called the multiplicities of 
$\alpha_1,\ldots, \alpha_r.$
      
2)  A quasi-parabolic structure on $W$ at $p_1,\ldots , p_n$ is just the 
condition a) in Definition \ref{def1} above at each $p=p_i.$ 
\erem

\bdefe\label{def2}Let $W_1$ and $W_2$ be two parabolic vector bundles on 
$X$ with parabolic structure at $p_1,\ldots , p_n.$
A morphism $G: W_1 \to W_2$ is a  vector bundle map from 
   $W_1$ to $W_2$ such that for any $p\in \{p_1,\ldots , p_n \},$ if we 
   denote by $g_p$ the linear map induced by $G$ on the fibers at $p,$ 
   we have 
   $$g_p(F_i(W_1)_p) \subset F_{j+1}(W_2)_p$$
    whenever $\alpha_i > \beta_j, $ where $\alpha_i$
   (resp.$\beta_j$ ) weights of $W_1$ (resp. $W_2.$)
\edefe

\bdefe\label{def3} 
Let $W$ be a parabolic vector bundle on $X$ with 
parabolic weights $\alpha_{1,i},\ldots, \alpha_{r_i,i}$  
with multiplicities $k_{1,i},\ldots, k_{r_i,i}$ at 
$p_i$ for $ i=1,\ldots,n.$  Then the parabolic degree of $W$
   is defined by
\[\text{Para deg} (W)=\text{deg}( W )+ \Sigma_i(\Sigma_j k_{j,i}\alpha_{j,i})\]
\edefe

\bdefe\label{def4}
Given a parabolic vector bundle $W$ and given any sub vector bundle $V$ and quotient bundle $V'$ 
then one defines in a natural way a parabolic structure on $V$ and $V'$
as follows:  for  $p \in \{p_1,\ldots, p_n \}$ the flag on $V_{p}$ (resp. $V'_{p}$) 
is induced by taking intersection (resp. quotient) with the flag of $W_{p_i}$ 
and weight attached for the subspace $F_k(V_{p})$ 
(resp. $F_k(V'_{p})$) is $\beta_k = \alpha_i$ where $i$ is the largest
(resp. largest) integer such that $F_k(V_{p})\subset F_i(W_{p})$
(resp. $F_i(W_{p})\to F_k(V_{p})$ is onto).
\edefe

\bdefe\label{def5}
   A parabolic vector bundle $V$ is said to be 
\em{Paraboic stable} (resp. \em{Parabolic semi-stable}) if 
for every parabolic sub-bundle $W$ of $V$ we have:
\[\frac{\text{Para deg} W}{\text{rank}(W) } < \frac{\text{Para deg} V}{\text{rank}(V)} \,\,(\text{resp.} \leq ). \]
\edefe

For details we refer to \cite{MS}
\section{Narasimhan-Seshadri type result for Paraboic bundles }

\bth\label{thm1} {\em Mehta-Seshadri Theorem}:
(1) Let $X$ be a smooth projective curve over 
an algebraically closed field with $g(X)\geq 2.$ 
Let $S$ be the set of all parabolic semi-stable bundles of rank $k,$ 
with fixed parabolic structure at a given point $p \in X,$ 
fixed weights $0 < \alpha_1<\ldots < \alpha_r<1$ of fixed multiplicities, 
fixed degree d and parabolic degree $0.$ 
Assume $\alpha_i$ are all rational. Two bundles $W,W′\in S$ are 
termed "equivalent'' if $\text{gr}W=\text{gr}W′.$ 
Then the set of equivalence classes of $S$ carries in a natural way 
the structure of a normal projective 
variety of dimension $k^2(g−1)+1+\text{dim}F$ where $F$ is the 
flag variety determined by the 
multiplicities of the parabolic structure at $p.$

(2) Assume that the field is field of complex numbers and $\Gamma$ be a
 discrete subgroup of $\text{PSL}(2, \mathbb{R})$ with a single equivalence class of cusp in
$\mathbb{R}\cup\infty$
for its action on the upper half plane $H.$ Let 
$$H^+ = H\cup \{Q \in \mathbb{R}| Q \text{is a 
cusp for} \,\,\Gamma \} $$
Then $X = H^+/{\Gamma}$  is a compact Riemann Surface
 Fix a cusp $Q \in H^+$ and  $P \in X$ be the corresponding point. Let $\Gamma_Q$ be the stabilizer. 
 Then the above variety is isomorphic to the set of equivalence classes of unitary representations 
 of $\Gamma$ with the image of the generator of $\Gamma_Q$ being conjugate to the diagonal matrix 
$({\exp}(2\pi i\alpha_1),\ldots,{\exp}(2\pi i\alpha_r))$ where each $\alpha_i$ is repeated $k_i$ times. 
 In particular, irreducible representations correspond to parabolic stable bundles.
\eeth

{\bf Applications of Mehta-Seshadri Theorem}
 
 The notion of parabolic structure has many applications in algebraic geometry
 For example they appear in  Hecke correspondences of Narasimhan and Ramanan
 \cite{NR}. It also, appear in the work of Nagaraj and  Seshadri on 
 moduli of torsion sheaves on curve which is a union of two smooth curves meeting at point \cite{NS}. 
 Parabolic bundles have had interesting applications 
 even outside algebraic geometry, in topology and physics. 
 One can find more details about the application of parabolic bundles by searching in the "web".

{\bf Generalizations}

1) {\em Parabolic vector bundles on Higher dimension varieties.}
  
Maruyama and Yokagava generalize the notion of parabolic bundles and 
 various associated notions from curves to non singular projective varieties of arbitrary dimension \cite{MY}. 
 The parabolic data now reside over an effective divisors. They proceed to construct a coarse 
 moduli scheme for stable parabolic sheaves over a non singular projective variety.  
 
2) {\em Generalized parabolic vector bundles.}
    
Usha  Bhosle, Narasimhan and Ramadas (see \cite{Uh1} \cite{NRa})
generalized the notion of parabolic bundles to 
what are called Generalized parabolic bundles and constructed moduli 
space of such bundles.
Generalized parabolic bundles are useful in the study of torsion free 
sheaves on nodal curves.
This notion is used by Usha Bhosle, Narasimhan and Ramadas,  Sun Xiaotao, 
Nagaraj and Seshadri, Ivan Kausz and many others in their work (see 
\cite{Ka}, \cite{NS}, \cite{NRa},\cite{Uh1}, \cite{Uh2}, \cite{Su}). 

3) {\em Parabolic Principal  bundles }

Several attempts were made to generalize the concept of Parabolic
bundles to $G$ bundles, where $G$ is a reductive algebraic group. 
Laszlo and Sorger defined the notion of Parabolic $G$ bundles and 
studied the Picard group of
moduli stack of such bundles on a projective smooth curve (see \cite{LS}). 
Balaji, Biswas and Nagaraj defined
   Parabolic $G$ bundle as a functor from catagory of $G$ modules to catagory 
   of Parabolic Vector
   bundles satisfying some conditions (following Nori's approch to principal 
   $G$ bundles) and
   latter as a ramified geometric objects (see, \cite{BBN1} \cite{BBN2}).

Let $ X$ be an irreducible smooth projective algebraic curve of 
   genus $ g \geq 2$ over the ground field 
  $ \mathbb{C}$, and let $ G$ be a semisimple simply connected algebraic group. 
  Balaji and Seshadri  introduce the notion of semistable and stable  parahoric torsors 
  under a certain Bruhat-Tits parahoric group scheme $ \mathcal G$ over $X.$ They construct the moduli space of 
  semistable parahoric $ \mathcal G$-torsors; 
  and identified the underlying topological space of this moduli space with 
  certain spaces of homomorphisms of Fuchsian groups into a maximal compact subgroup of $ G$. 
  The results give a generalization of the earlier results of Mehta and Seshadri on parabolic vector bundles(See \cite{BS}). 
  
\section{Parabolic Higgs bundles}

{\bf Motivation}:
Ordinary Higgs bundles on a Riemann Surface were introduced by Hitchin. Recall, if $X$
is a compact Riemann, then Higgs bundle on $X$ is a pair $(V, \theta)$ where $V$ is a vector bundle on $X$ and 
$\theta : V \to V\otimes\Omega^1_X$ a homomorphism of bundles, where $\Omega^1_X$ is the bundle of holomorphic
1-forms on $X.$ There is a notion (semi-) stability associated these pairs which is natural generalization of
same notions for a vector bundle. Moduli space of Higgs bundles exists as quasi projective variety.
There is an analogous Narasimhan-Seshadri theorem which sets up correspondence
between certain Higgs bundles on $X$ and linear representations of the fundamental group of $X.$ In order to generalize
Mehta-Seshadri correspondence one need to extend the definition of Higgs bundle to the Parabolic case.

{\bf Parabolic Higgs bundles}
 
 Paraboic Higgs bundle on compact connected Riemann $X$ is a pair
 $(V, \theta)$ where $V$ is a parabolic vector bundle on $X$ and $\theta $ is 
 as indicated in the $\theta : V \to V\otimes\Omega^1_X$ a homomorphism of 
 bundles which at each parabolic point preserves the flag.
 
One can define (semi-)stability of Parabolic Higgs bundles and moduli of such
 objects were constructed by Konno Hiroshi, Maruyama-Yokogawa and many
 others. Analog of Mehta-Seshadi
 theorem holds in this context (see \cite{Ko}, \cite{MY} \cite{Yo}).
  
In \cite{Ni} Nitsure generalized the notion of Higgs bundles to what are now
called Hitchin Pairs and constructed the moduli spaces stable Hitchin Pairs.
This  notion of Hitchin pairs has been generalized to  parabolic Hitchin pairs.
The moduli spaces of Parabolic Higgs bundles is still a active area of 
research.


\begin{thebibliography}{9999}
\bibitem{BBN1} Balaji, Vikraman; Biswas, Indranil and Nagaraj, D. S.
\textit{Principal bundles over projective manifolds with parabolic structure 
over a divisor.} Tohoku Math. J. (2) 53 (2001), no. 3, 337–367.

\bibitem{BBN2} Balaji, V; Biswas, I and Nagaraj, D. S.\textit{Ramified G-
bundles as parabolic bundles.} J. Ramanujan Math. Soc. 18 (2003), 
no. 2, 123–138.

\bibitem{BS} Balaji, V.; Seshadri, C. S. 
\textit{Moduli of parahoric 𝒢-torsors on a 
compact Riemann surface.} J. Algebraic Geom. 24 (2015), no. 1, 1–49.

\bibitem{Hi} Hitchin, N. J. \textit{The self-duality equations on a 
Riemann surface.} Proc. London Math. Soc. (3) 55 (1987), no. 1, 59–126. 

\bibitem{Ka}  Kausz, Ivan. \textit{A Gieseker type degeneration of moduli 
stacks vector bundles on curves.} Trans. Amer. Math. Soc. 357 (2005), 
no. 12, 4897–4955.

\bibitem{Ko} Konno, Hiroshi \textit{Construction of the moduli space of 
stable parabolic Higgs bundles on a Riemann surface.}
J. Math. Soc. Japan 45 (1993), no. 2, 253–276. 

\bibitem{LS} Y. Laszlo and C. Sorger, \textit{The line bundles on the moduli 
of parabolic G-bundles over curves and their sections.} Ann. Sci. École Norm. 
Sup. 30 (1997), 499–525 

\bibitem{MS} V.B. Mehta and C.S. Seshadri, \textit{Moduli of vector
bundles on curves with Parabolic Structure} Math. Ann. 248, 205-239 (1980)

\bibitem{MY} M. Maruyama and K. Yokogawa, \textit{Moduli of parabolic stable 
sheaves}, Math. Ann. 293 (1992), 77--99.

\bibitem{NS} D.S. Nagaraj and C.S. Seshadri, \textit{Degenerations of the
 moduli spaces of vector bundles on curves I} Proc. Indian Aced.Sci.
 (Math.Sci). Vol. 107, No.2, May 1997,pp.101-137.

\bibitem{NRa}  Narasimhan, M. S.; Ramadas, T. R. \textit{Factorisation of 
generalised theta functions. I.} Invent. Math. 114 (1993), no. 3, 565–623. 
 
\bibitem{NR} Narasimhan, M. S. and  Ramanan, S. \textit{Geometry of Hecke
 cycles. I.} C. P. Ramanujam—a tribute, pp. 291–345, Tata Inst. Fund. 
 Res.Studies in Math., 8, Springer, Berlin-New York, 1978. 

\bibitem{Ni} Nitsure, Nitin. \textit{Moduli space of semistable 
pairs on a curve.}
Proc. London Math. Soc. (3) 62 (1991), no. 2, 275–300. 

\bibitem{Uh1} Bhosle, Usha. \textit{Generalised parabolic bundles and 
applications to torsionfree sheaves on nodal curves.} Ark. Mat. 30 (1992), 
no. 2, 187–215.
 
\bibitem{Uh2}  Bhosle, Usha. \textit{Generalized parabolic sheaves on an 
integral projective curve.} Proc. Indian Acad. Sci. Math. Sci. 102 (1992), 
no. 1, 13–22. 

\bibitem{Su}  Sun, Xiaotao. \textit{Degeneration of moduli spaces and generalized theta functions.} J. Algebraic Geom. 9 (2000), no. 3, 459–527. 

\bibitem{Yo} Yokogawa, Kôji. \textit{Moduli of stable pairs.}
J. Math. Kyoto Univ. 31 (1991), no. 1, 311–327. 

\end{thebibliography}
\end{document}